\documentclass[11pt]{article}
\usepackage{amsmath}
\usepackage{amssymb}
\usepackage{latexsym}
\usepackage{amsthm}
\usepackage{tabularx}
\usepackage{booktabs}
\usepackage{mathrsfs}
\usepackage{enumitem}
\usepackage{ytableau}
\usepackage{graphicx}
\usepackage{algorithm,algorithmic}
\usepackage[colorlinks,
            linkcolor=red,
            anchorcolor=blue,
            citecolor=green]{hyperref}

\newtheorem{thm}{Theorem}[section]
\newtheorem{lem}[thm]{Lemma}
\newtheorem{prop}[thm]{Proposition}

\allowdisplaybreaks

\begin{document}
\begin{center}
{\large \bf   A bijective proof of the hook-length formula \\
for standard immaculate tableaux}
\end{center}

\begin{center}
Emma L.L. Gao$^1$ and Arthur L.B. Yang$^{2}$\\[6pt]

Center for Combinatorics, LPMC-TJKLC\\
Nankai University, Tianjin 300071, P. R. China\\[6pt]

Email: $^{1}${\tt gaolulublue@mail.nankai.edu.cn},
       $^{2}${\tt yang@nankai.edu.cn}
\end{center}

\noindent\textbf{Abstract.}
In this paper, we present a direct bijective proof of the hook-length formula for standard immaculate tableaux, which arose in the study of  non-commutative symmetric functions.
Our proof is along the spirit of Novelli, Pak and Stoyanovskii's combinatorial proof of the hook-length formula for standard Young tableaux.

\noindent \emph{AMS Classification 2010:} 05E05


\section{Introduction}

In the study of non-commutative symmetric functions, Berg, Bergeron, Saliola, Serrano and Zabrocki \cite{BBSSZ} introduced the notion of immaculate tableaux, which was indexed by compositions of integers.
They obtained an amazingly simple product formula to enumerate standard immaculate tableaux, which is analogous to the hook-length formula for standard Young tableaux. Their proof is by induction on the length of the composition. The objective of this paper is to give a direct bijective proof of the hook-length formula for standard immaculate tableaux.

The classical hook length formula for standard Young tableaux was first discovered by Frame, Robinson and Thrall \cite{FRT}. To explore why hooks appear in this formula, many proofs have been published based on different methods. The first step towards this direction was given by Hillman and Grassl \cite{HG}. They proved a special case of Stanley's hook-content formula, from which the hook-length formula follows. Later, Greene, Nijenhuis and Wilf \cite{GNW} found a probabilistic proof using the hook walk which shows clearly the role of hooks. The first bijective proof was given by Franzblau and Zeilberger \cite{FZ}, though it is not so direct. Even though there were so many proofs, none of them was considered satisfactory. Novelli, Pak and Stoyanovskii \cite{NPS} presented an elegant bijective proof of the hook-length formula, based on the work of Pak and Stoyanovskii \cite{PS}.

Motivated by Novelli, Pak and Stoyanovskii's combinatorial proof of the classical hook-length formula, it is natural to consider whether a naturally bijective proof exists for the hook-length formula of standard immaculate tableaux, which could clearly illuminate the role of hooks. In this paper, we shall present such a proof.

Let us first review some notation and terminology concerning the hook-length formula for standard immaculate tableaux. A \emph{composition} $\alpha$ of a positive integer $n$, denoted by $\alpha\models n$,  is a tuple $\alpha=(\alpha_1,\alpha_2,\ldots,\alpha_k)$ of positive integers such that $\sum_{i=1}^k\alpha_i=n$. The entries $\alpha_i$ are called the \emph{parts} of $\alpha$, and the number of parts is called the \emph{length} of $\alpha$, denoted by $\ell(\alpha)$.
Each composition is associated to a diagram of left-justified array of cells. Given a composition $\alpha=(\alpha_1,\alpha_2,\ldots,\alpha_k)$,  the corresponding diagram has $\alpha_i$ cells in the $i$-th row. Here we number the rows from top to bottom and the columns from left to
right. The cell in the $i$-th row and $j$-th column is denoted by the pair $(i, j)$.
For example, the diagram of the composition $(4,1,2,3)$ is as follows.

\begin{figure}[h!]
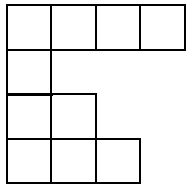

\centering
 \ydiagram{4,1,2,3}
\caption{The diagram of the composition $(4,1,2,3)$}
\label{fig:  diagram example}
\end{figure}

Following Berg, Bergeron, Saliola, Serrano and Zabrocki \cite{BBSSZ}, we now introduce the definitions of hooks and immaculate tableaux. Given a composition $\alpha$ and a cell $c=(i,j)$ in $\alpha$, the \emph{hook} of $c$, denoted $H_c$, is defined by
\begin{align*}
H_c =H_{i,j}= \left\{
             \begin{array}{ll}
             {\{(i',j'):i\leq i'\leq \ell(\alpha), 1\leq j'\leq \alpha_{i'}\}}, &\mbox{if $j=1$;} \\[5pt]
             {\{(i,j'):j\leq j'\leq \alpha_i\}}, &\mbox{if $j>1$.}
             \end{array}
        \right.
\end{align*}
Correspondingly, the \emph{hook{-}length} of the cell $c=(i,j)$, denoted by $h_c$, is defined as
$$h_c=h_{i,j}=|H_{i,j}|.$$
For example, taking the cells $(1,2)$ and $(2,1)$ of the composition $\alpha=(4,1,2,3)$, the hooks $H_{1,2}$ and $H_{2,1}$ are depicted in Figure \ref{fig:  hook example} as the sets of dotted cells. Clearly, we have $h_{1,2}=3$ and $h_{2,1}=6$.

\begin{center}
\makeatletter\def\@captype{figure}\makeatother
\begin{minipage}{.3\textwidth}
\centering
\ydiagram{4,1,2,3}
*[\bullet]{0,1+3}\\[6pt]
$H_{1,2}$
\end{minipage}
\makeatletter\def\@captype{figure}\makeatother
\begin{minipage}{.3\textwidth}
\centering
\ytableausetup{centertableaux}
\ytableaushort
{\none,\bullet, \bullet \bullet,\bullet \bullet \bullet}
* {4,1,2,3}
* {4} \\[6pt]
$H_{2,1}$
\end{minipage}
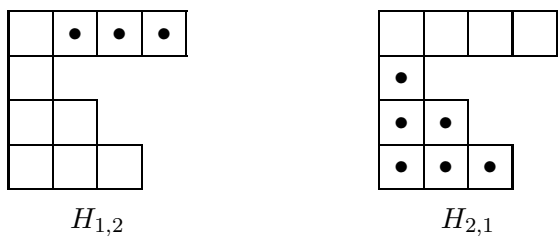
\caption{The hooks of $(1,2)$ and $(2,1)$}
\label{fig:  hook example}
\end{center}

We proceed to introduce the concept of immaculate tableaux. Given a composition $\alpha \models n$,
a \emph{tableau} of shape $\alpha$ is defined to be an array $T=(T_{ij})$ obtained by filling the diagram of $\alpha$ with positive integers. For the convenience, let $T_{ij}=\infty$ if the cell $(i,j)$ does not belong to the diagram of $\alpha$.
We say that $T$ has content $\beta=(\beta_1,\beta_2,\ldots)$ if, for each $i\geq 1$, there are $\beta_i$ entries equal to $i$ in the array. For any $1\leq i\leq \ell(\alpha)$ and $1\leq j\leq \alpha_i$, the entry $T_{ij}$ is said to be {\emph{stable}} in $T$ if either of the following two conditions holds:
\begin{itemize}
\item $j>1$ and $T_{ij}\leq T_{i,j+1}$; or
\item $j=1$ and $T_{ij}<T_{i+1,j}$ and $T_{ij}\leq T_{i,j+1}$.
\end{itemize}
An \emph{immaculate tableau} of shape $\alpha$ is a tableau $T=(T_{ij})$ of shape $\alpha$ such that:
\begin{itemize}
\item[(i)] all entries $T_{ij}$ are stable; and

\item[(ii)] for any $m\geq 2$, if $m$ appears in $T$, so does $m-1$.
\end{itemize}
Condition (ii) implies that the content of an immaculate tableau must be a composition $\beta$ of the form
$(\beta_1,\ldots,\beta_{\ell(\beta)})$ whose components are all positive.
Given an immaculate tableau $T$ of shape $\alpha$, we say that it is \emph{standard} if $T$ has content $(1^n)$.

Remarkably, Berg, Bergeron, Saliola, Serrano and Zabrocki \cite{BBSSZ} found that the standard immaculate tableaux can be enumerated by using the
above defined hook-lengths of the cells of the indexed composition.
Given a composition $\alpha$, denote by $f^\alpha$ the number of standard immaculate tableaux of shape $\alpha$. The hook-length formula for standard immaculate tableaux reads as follows.

\begin{thm}[{\cite[Proposition 3.13]{BBSSZ}}] \label{hook_immaculate}
For any composition $\alpha\models n$, we have
\begin{align}\label{hook_formula}
 f^\alpha=\frac{n!}{\prod_{(i,j)\in \alpha}h_{i,j}}.
 \end{align}
\end{thm}

For instance, the composition $(2, 1, 2)$ has hook lengths given by

\begin{center}
\begin{ytableau}
  5& 1 \\
  3\\
 2& 1
\end{ytableau}
\end{center}
From the hook-length formula \eqref{hook_formula} it follows that the number of standard immaculate tableaux of shape $(2, 1, 2)$ is
$$\frac{5!}{5\cdot1\cdot3\cdot2\cdot1}=4.$$
In fact, there are exactly $4$ standard immaculate tableaux of shape $(2, 1, 2)$ as illustrated below.

\begin{center}
\makeatletter\def\@captype{figure}\makeatother
\centering
\begin{minipage}{.15\textwidth}
\begin{ytableau}
  1& 2  \\
  3\\
 4& 5
\end{ytableau}
\end{minipage}
\makeatletter\def\@captype{figure}\makeatother
\begin{minipage}{.15\textwidth}
\begin{ytableau}
  1& 3  \\
  2\\
 4& 5
\end{ytableau}
\end{minipage}
\makeatletter\def\@captype{figure}\makeatother
\begin{minipage}{.15\textwidth}
\begin{ytableau}
  1& 4  \\
  2\\
 3& 5
\end{ytableau}
\end{minipage}
\makeatletter\def\@captype{figure}\makeatother
\begin{minipage}{.15\textwidth}
\begin{ytableau}
  1& 5  \\
  2\\
 3& 4
\end{ytableau}
\end{minipage}
\caption{The standard immaculate tableaux of shape $(2, 1, 2)$}
\label{fig: standard immaculate tableaux}
\end{center}

\section{A bijective proof of Theorem
\ref{hook_immaculate}}

The aim of this section is to present  a bijective proof of Theorem
\ref{hook_immaculate}. Towards this end, we first rewrite (\ref{hook_formula}) as
\begin{align}\label{eq-rewritten}
n!=f^\alpha\prod_{(i,j)\in \alpha}h_{i,j}.
\end{align}
Then we need to construct two sets such that their cardinalities are respectively given by the left-hand side and the right-hand side of \eqref{eq-rewritten}.
Let $X$ be the set of  tableaux of shape $\alpha$ and content $(1^n)$, and let $Y$ be the set of $\{(P, J)\}$, where $P$ is a standard immaculate tableau of shape $\alpha$ and $J$ is an array of shape $\alpha$ with $J_{i,j}\in \{1,\ldots, h_{i,j}\}$. We call $J$ a hook tableau. It is easy to see that
\begin{align*}
|X|=n!\qquad \mbox{  and  } \qquad |Y|=f^\alpha\prod_{(i,j)\in \alpha}h_{i,j}.
\end{align*}
There remains to show that there exists a bijection between $X$ and $Y$.

For the construction of our bijection, a total order on the cells of $\alpha$ is needed. Following
Novelli, Pak and  Stoyanovskii, we totally order the cells of the diagram of $\alpha$
by reverse lexicographic order on their coordinates.
Precisely, we have
\begin{center}
  $(i,j)\leq (i',j')$ if and only if $j>j'$; or $j=j'$ and $i\geq i'.$
\end{center}
As will be shown below, this order is critical for the construction of our bijection.
Label the cells of $\alpha$ in the given order
$c_1<c_2<\cdots<c_n.$
For example, Figure \ref{The total order} displays the total order for the diagram of shape $(4,1,4,2,1)$.

\makeatletter\def\@captype{figure}\makeatother
\begin{minipage}{.99\textwidth}
\centering
\begin{ytableau}
c_{12}& c_7& c_4 & c_2\\
c_{11}\\
c_{10}&  c_6 & c_3 & c_1\\
c_9&c_5 \\
c_8
\end{ytableau}
\end{minipage}
\caption{The total order}
\label{The total order}
Given a tableau $T$ of shape $\alpha$ and a cell $c$, let $T^{\leq c}$ (resp. $T^{<c}$) denote the partial tableau composed of all cells $b$ of $T$ with $b\leq c$ (resp. $b<c$). For example,
\begin{center}
\begin{align*}
\mbox{if \qquad} T = \begin{ytableau}
    11& 5& 8& 9\\
  3\\
 10&  2 & 4 & 12\\
 1&6 \\
 7
\end{ytableau},
\qquad \mbox{then \qquad}
T^{\leq c_8} = \begin{ytableau}
  \none& 8& 5 & 9\\
   \none\\
  \none&  12 & 2 & 4\\
  \none&6 \\
 7
\end{ytableau}.
\end{align*}
\end{center}
For the convenience, we say that the partial tableau $T^{\leq c}$ (resp. $T^{<c}$) is \emph{standard} if all the entries in $T^{\leq c}$ (resp. $T^{<c}$) are stable with respect to the diagram of $\alpha$.

We now construct a map $\psi$ from $Y$ to $X$.
Given a pair $(P,J)\in Y$, we construct a tableau $T\in X$ in the following way.
Without loss of generality, we may assume that $n>1$. Begin with $(P_1,J_1)=(P,J)$.
If $(P_k,J_k)$ are defined for $1\leq k<n$, then let $J_{k+1}=J_k$ except for $(J_{k+1})_{ij}=1$ if the cell $c_{n+1-k}$ lies in the $i$-th row and $j$-th column. Suppose that the $(J_{k})_{ij}$-th cell of the hook set $H_{ij}$, reading from left to right and top to bottom, lies in the $i'$-th row and $j'$-th column of the diagram of $\alpha$. The cells $(i,j)$ and $(i',j')$ uniquely determine a path $L$ in the following way.
\begin{itemize}
\item[(a)] If $i=i'$, then let $L=\{(i,j), (i,j+1),\ldots,(i,j')\}$;

\item[(b)] Suppose that $i\neq i'$. Then by the definition of the hook function, we must have $i<i'$ and $j=1$. In this case, let $L=\{(i,1), (i+1,1),\ldots,(i',1),(i',2),\ldots,(i',j')\}.$
\end{itemize}
The tableau $P_{k+1}$ is obtained from $P_k$ by a circular right shift of the entries on the path $L$.
If the process ends at $(P_n,J_n)$, then all entries of $J_n$ are $1$, and let $T=P_n$. The map $\psi$ is defined by $\psi(P,J)=T$. Note that $P_1^{\leq c_n}=P$ is standard. A moment's thought shows that the partial tableau $P_k^{\leq c_{n+1-k}}$ is standard for any $k$.

For example, for $\alpha=(4,1,4,2,1)$, let
\begin{center}
\begin{tabular}{@{}llc@{}}
P = \begin{ytableau}
    1& 5& 8& 9\\
  2\\
 3&  4 & 11& 12\\
 6&10 \\
 7
\end{ytableau}&
and &
J = \begin{ytableau}
    8& 2& 1& 1\\
  3\\
6&  3 & 1 &1\\
 1&1 \\
 1
\end{ytableau}.
\end{tabular}
\end{center}
We have a sequence of pairs $\{(P_i, J_i)\}_{i=1}^{12}$ as follows, where the entries in the path $L$ are underlined at each step.

\ytableausetup{mathmode, boxsize=2em, smalltableaux}
\begin{center}
\begin{tabular}{@{}ccc}
%\toprule
$i$ &$P_i$ & $J_i$\\
\midrule
 $1$&\begin{ytableau}
    \underline{1}& 5& 8& 9\\
  \underline{2}\\
 \underline{3}& \underline{ 4} & \underline{11}& 12\\
 6&10 \\
 7
\end{ytableau} & \begin{ytableau}
    8& 2& 1& 1\\
  3\\
6&  3 & 1 &1\\
 1&1 \\
 1
\end{ytableau}\\
\\
$2$&\begin{ytableau}
    11& 5& 8& 9\\
  \underline{1}\\
 \underline{2}& \underline{ 3} & 4 & 12\\
 6&10 \\
 7
\end{ytableau} & \begin{ytableau}
    1& 2& 1& 1\\
  3\\
6&  3 & 1 &1\\
 1&1 \\
 1
\end{ytableau}\\
\\
$3$&\begin{ytableau}
    11& 5& 8& 9\\
  3\\
 \underline{1}&  2 & 4 & 12\\
 \underline{6}&\underline{10} \\
 7
\end{ytableau} & \begin{ytableau}
    1& 2& 1& 1\\
  1\\
6&  3 & 1 &1\\
 1&1 \\
 1
\end{ytableau}\\
\\
$4$&\begin{ytableau}
    11& 5& 8& 9\\
  3\\
 10&  2 & 4 & 12\\
 \underline{1}&6 \\
 7
\end{ytableau} & \begin{ytableau}
    1& 2& 1& 1\\
  1\\
 1&  3 & 1 &1\\
 1&1 \\
 1
\end{ytableau}
\end{tabular}
\quad $\Rightarrow$ \quad
\begin{tabular}{@{}ccc}
%\toprule
$i$ &$P_i$ & $J_i$\\
\midrule
$5$&\begin{ytableau}
    11& 5& 8& 9\\
  3\\
 10&  2 & 4 & 12\\
 1&6 \\
 \underline{7}
\end{ytableau} & \begin{ytableau}
    1& 2& 1& 1\\
  1\\
 1&  3 & 1 &1\\
 1&1 \\
 1
\end{ytableau}\\
\\ $6$&\begin{ytableau}
    11& \underline{5}& \underline{8}& 9\\
  3\\
 10&  2 & 4 & 12\\
 1&6 \\
 7
\end{ytableau} & \begin{ytableau}
    1& 2& 1& 1\\
  1\\
 1&  3 & 1 &1\\
 1&1 \\
 1
\end{ytableau}\\
\\ $7$&\begin{ytableau}
    11& 8& 5 & 9\\
  3\\
 10&  \underline{2} & \underline{4} & \underline{12}\\
 1&6 \\
 7
\end{ytableau} & \begin{ytableau}
    1& 1& 1& 1\\
  1\\
 1&  3 & 1 &1\\
 1&1 \\
 1
\end{ytableau}\\
\\ $8$ &
\begin{ytableau}
    11& 8& 5 & 9\\
  3\\
 10&  12 & 2 & 4\\
 1&\underline{6} \\
 7
\end{ytableau} & \begin{ytableau}
    1& 1& 1& 1\\
  1\\
 1&  1 & 1 &1\\
 1&1 \\
 1
\end{ytableau}
\end{tabular}
\quad $\Rightarrow$ \quad
\begin{tabular}{@{}ccc}
%\toprule
$i$ &$P_i$ & $J_i$\\
\midrule
$9$ &
\begin{ytableau}
    11& 8& \underline{5} & 9\\
  3\\
 10&  12 & 2 & 4\\
 1&6 \\
 7
\end{ytableau} & \begin{ytableau}
    1& 1& 1& 1\\
  1\\
 1&  1 & 1 &1\\
 1&1 \\
 1
\end{ytableau}\\
\\ $10$ &
\begin{ytableau}
    11& 8& 5 & 9\\
  3\\
 10&  12 & \underline{2} & 4\\
 1&6 \\
 7
\end{ytableau} & \begin{ytableau}
    1& 1& 1& 1\\
  1\\
 1&  1 & 1 &1\\
 1&1 \\
 1
\end{ytableau}\\
\\ $11$ &
\begin{ytableau}
    11& 8& 5 & \underline{9}\\
  3\\
 10&  12 & 2 & 4\\
 1&6 \\
 7
\end{ytableau} & \begin{ytableau}
    1& 1& 1& 1\\
  1\\
 1&  1 & 1 &1\\
 1&1 \\
 1
\end{ytableau}\\
\\ $12$ &
\begin{ytableau}
    11& 8& 5 & 9\\
  3\\
 10&  12 & 2 & 4\\
 1&6 \\
 7
\end{ytableau} & \begin{ytableau}
    1& 1& 1& 1\\
  1\\
 1&  1 & 1 &1\\
 1&1 \\
 1
\end{ytableau}%\bottomrule
\end{tabular}
\end{center}

\ytableausetup{mathmode, boxsize=normal}

Our main result is as follows.

\begin{thm} \label{thm-main} The map $\psi$ is a bijection from $Y$ to $X$.
\end{thm}

To show that $\psi$ is a bijection, it suffices to construct a map $\phi$ from $X$ to $Y$ such that $\phi=\psi^{-1}$. This map $\phi$ is based on a modified jeu de taquin performed on $X$. Suppose that we are given a tableau $T$ of shape $\alpha$ in $X$. For any $1\leq e\leq n$, denote by $(i,j)$ the unique cell in $T$ such that $T_{ij}=e$. To each $e$, we will associate a transformation $\mathrm{jdt}_e(T)$ of $T$ called a modified jeu de taquin slide of $T$ with respect to $e$. If $e$ is stable, then we do nothing. If $e$ is not stable, then there are two cases to consider according to whether $e$ lies in the first column of the diagram of $\alpha$.
\begin{itemize}
\item[(a)] If $T_{ij}=e$ for some $1<j\leq\alpha_i$, then interchange $T_{ij}$ and $T_{i,j+1}$;

\item[(b)] If $T_{i1}=e$ for some $1\leq i\leq \ell(\alpha)$, then interchange $T_{ij}$ and the smaller one of $\{T_{i,j+1}, T_{i+1,j}\}$.
\end{itemize}
Then look at the stability of $e$, and repeat the same procedure. This will eventually terminate since the entry $e$, if it moves, will move either downwards or rightwards in $T$ at each step. The path of $e$ in $T$ is defined as the set of the cells that $e$ passes through when applying $\mathrm{jdt}_e(T)$. By convention, we also denote the resulting tableau by $\mathrm{jdt}_e(T)$.

We use an example to illustrate the modified jeu de taquin algorithm. Taking $e=10$ and
\begin{align*}
T = \begin{ytableau}
    11& 5& 8& 9\\
  3\\
 10&  2 & 4 & 12\\
 1&6 \\
 7
\end{ytableau},
\end{align*}
then the process of $\mathrm{jdt}_e(T)$ is as follows,
\[
\begin{matrix}
T & = &
& \begin{ytableau}
    11& 5& 8& 9\\
  3\\
 \textbf{10}&  2 & 4 & 12\\
 \underline{1}&6 \\
 7
\end{ytableau}
& \rightarrow &
%\raisebox{-5pt}
{\begin{ytableau}
    11& 5& 8& 9\\
  3\\
 1&  2 & 4 & 12\\
 \textbf{10}&\underline{6} \\
 7
\end{ytableau}}
 &\rightarrow &
%\raisebox{-15pt}
{\begin{ytableau}
    11& 5& 8& 9\\
  3\\
 1&  2 & 4 & 12\\
 6&\textbf{10} \\
 7
\end{ytableau} }
& = & \mathrm{jdt}_e(T),
\end{matrix}
\]
where the entry $e$ is in boldface and the integers interchanged with $e$ are underlined. Clearly, the path of $10$ in $T$ is
$\{(3,1),(4,1),(4,2)\}$.

The following result is evident, and we omit the straightforward details.

\begin{lem}\label{lemma-path}
Given $T\in X$ and $1\leq e\leq n$, let $(i,j)$ be the unique cell in $T$ such that $T_{ij}=e$. If $j>1$, then the path of $e$ in $T$ is of the form
\begin{align*}%\label{eq-1}
\{(i,j),(i, j+1),\ldots,(i,j+k)\}, \quad \mbox{for some $k\geq 0$.}
\end{align*}
If $j=1$, then the path of $e$ in $T$ is of the form
\begin{align*}%\label{eq-2}
\{(i,1), (i+1, 1),\ldots,(i+k,1),(i+k, 2),\ldots,(i+k,l)\}, \quad \mbox{for some $k\geq 0$ and $l\geq 1$.}
\end{align*}
Moreover, the tableau $\mathrm{jdt}_e(T)$ is obtained from $T$ by a circular left shift of the entries on the
path of $e$.
\end{lem}

The next result shows that the modified jeu de taquin slide preserves the stability of most entries in the tableau.

\begin{prop}\label{standard pro}
Given $T\in X$ and $1\leq e\leq n$, denote by $c$ the unique cell in $T$ such that $T_{c}=e$. Then if $T^{<c}$ is standard, so is $\mathrm{jdt}_e(T)^{\leq c}$.
\end{prop}

\proof Note that $\mathrm{jdt}_e(T)^{\leq c}=T^{\leq c}$ except for the entries of the cells in the path of $e$ in $T$.  Suppose that $c$ is the $(i,j)$ cell of $T$. There are two cases to consider.

\begin{itemize}
\item[(a)] If $j>1$, then $\mathrm{jdt}_e(T)^{\leq c}$ has no cell in the first column of the diagram of $\alpha$. We only need to show that the entries of $\mathrm{jdt}_e(T)^{\leq c}$ are strictly increasing in each row. By Lemma \ref{lemma-path}, the path of $e$ is of the form
    $$\{(i,j),(i, j+1),\ldots,(i,j+k)\}.$$
    Therefore, we have
    $$
     (\mathrm{jdt}_e(T))_{rs}=\left\{\begin{array}{ll}
    T_{r,s+1}, & \mbox{for $r=i$ and $j\leq s\leq j+k-1$},\\[3pt]
    T_{i,j}, & \mbox{for $r=i$ and $s=j+k$},\\[3pt]
    T_{rs}, & \mbox{otherwise.}
    \end{array}\right.
    $$
The entries of $\mathrm{jdt}_e(T)^{\leq c}$ in each row other than the $i$-th row, are identically the same as those of $T^{<c}$, and hence are increasing since $T^{<c}$ is standard. While for the $i$-th row, the form of the path of $e$ implies that
$$T_{i,j+1}<T_{i,j+2}<\cdots<T_{i,j+k}<T_{i,j}<T_{i,j+k+1},$$
that is
$$\mathrm{jdt}_e(T)_{i,j}<\mathrm{jdt}_e(T)_{i,j+1}<\cdots<\mathrm{jdt}_e(T)_{i,j+k-1}<\mathrm{jdt}_e(T)_{i,j+k}<T_{i,j+k+1}.$$
Thus, the entries in the $i$-th row of $\mathrm{jdt}_e(T)^{\leq c}$ are also increasing, as desired.

\item[(b)] If $j=1$, then $\mathrm{jdt}_e(T)^{\leq c}$ must contain a cell in the first column of the diagram of $\alpha$. By Lemma \ref{lemma-path},  the path of $e$ is of the form
\begin{align*}
\{(i,1), (i+1, 1),\ldots,(i+k,1),(i+k, 2),\ldots,(i+k,l)\}.
\end{align*}
Our proof will be divided into two subcases:

\begin{itemize}
\item[(b1)] The case of $l=1$. In this case, we have
    $$
    (\mathrm{jdt}_e(T))_{rs}=\left\{\begin{array}{ll}
    T_{r+1,1}, & \mbox{for  ${s}=1$ and $i\leq r\leq i+k-1$},\\[3pt]
    T_{i,1}, & \mbox{for ${s}=1$ and $r=i+k$},\\[3pt]
    T_{rs}, & \mbox{otherwise.}
    \end{array}\right.
    $$
Since $T^{<c}$ is standard, the form of the path of $e$ implies that
$$T_{i+1,1}<\cdots<T_{i+k,1}<T_{i,1}<T_{i+k+1,1}
<\cdots<T_{{\ell(\alpha)},1},$$
that is,
$$(\mathrm{jdt}_e(T))_{i,1}<\cdots<(\mathrm{jdt}_e(T))_{i+k,1}<(\mathrm{jdt}_e(T))_{i+k+1,1}<\cdots<(\mathrm{jdt}_e(T))_{{\ell(\alpha)},1}.$$
This means that the entries of $\mathrm{jdt}_e(T)^{\leq c}$ are strictly increasing in the first column.
There remains to show that each row of $\mathrm{jdt}_e(T)^{\leq c}$ are strictly increasing.
Since all entries of $T$ in any other than the first column remain fixed when applying $\mathrm{jdt}_e(T)$, it suffices to show that, for $i\leq r\leq i+k$,
$$(\mathrm{jdt}_e(T))_{r,1}<(\mathrm{jdt}_e(T))_{r,2}.$$
This is true, since, for $i\leq r\leq i+k-1$,
$$(\mathrm{jdt}_e(T))_{r,1}=T_{r+1,1}<T_{r,2}=(\mathrm{jdt}_e(T))_{r,2},$$
and
$$(\mathrm{jdt}_e(T))_{i+k,1}=T_{i,1}<T_{i+k,2}=(\mathrm{jdt}_e(T))_{i+k,2},$$
as implied by the form of the path of $e$. Therefore, the entries of $\mathrm{jdt}_e(T)^{\leq c}$ are increasing along each row.

\item[(b2)] The case of $l>1$.
In this case, we have
    $$
    (\mathrm{jdt}_e(T))_{rs}=\left\{\begin{array}{ll}
    T_{r+1,1}, & \mbox{for $s=1$ and $i\leq r\leq i+k-1$},\\[3pt]
    T_{i+k,2}, & \mbox{for $s=1$ and $r=i+k$},\\[3pt]
    T_{i+k,s-1}, & \mbox{for $1<s<l-1$ and $r=i+k$},\\[3pt]
    T_{i,1}, & \mbox{for $s=l$ and $r=i+k$},\\[3pt]
        T_{rs}, & \mbox{otherwise.}
    \end{array}\right.
    $$
    The strict increasing property of the $(i+k)$-th row of $\mathrm{jdt}_e(T)^{\leq c}$ can be proved by using similar arguments for case (a). While similar arguments for case (b1) could be used to show the strict increasing property of other rows of $\mathrm{jdt}_e(T)^{\leq c}$, as well as that of the first column of $\mathrm{jdt}_e(T)^{\leq c}$ except for the relation
    $\mathrm{jdt}_e(T)_{i+k,1}<\mathrm{jdt}_e(T)_{i+k+1,1}$. By the form of the path of $e$, we see that
    $$\mathrm{jdt}_e(T)_{i+k,1}=T_{i+k,2}<T_{i+k+1,1}=\mathrm{jdt}_e(T)_{i+k+1,1}.$$
    This completes the proof of the strict increasing property of the first column of $\mathrm{jdt}_e(T)^{\leq c}$.
\end{itemize}
\end{itemize}

Combining (a) and (b), we obtain the desired result.
\qed

We proceed to describe the inverse map $\phi: X\rightarrow Y$.  Suppose that $T\in X$ is a tableau of shape $\alpha$ and the cells of $\alpha$ are ordered as $c_1<c_2<\cdots<c_n$. We shall associate with $T$ a pair $(P,J)\in Y$, as follows. Begin with $(T_1,S_1)$, where $T_1=T$ and $S_1$ is the array of shape $\alpha$ with all entries equal to $1$. If $(T_{k}, S_{k})$ are defined for $1\leq k\leq n-1$, then let
$T_{k+1} =\mathrm{jdt}_{e_{k+1}}(T_{k})$, where $e_{k+1}$ is the entry of the cell $c_{k+1}$ in $T_{k}$. Suppose that the  path of $e_{k+1}$ in $T_{k}$ starts at $c_{k+1}=(i,j)$ and ends at $(i',j')$. Let $S_{k+1}=S_k$ except for the values
\begin{align}\label{eq-sk}
(S_{k+1})_{ij}= \left\{
             \begin{array}{lcl}
            j'-j+1, &\text{if} &i'=i; \\
            \alpha_i+\cdots+\alpha_{i'-1}+j', &\text{if} &i'>i.
             \end{array}
        \right.
\end{align}
Suppose that the process ends at $(T_{n}, S_{n})$. We claim that $(T_{n}, S_{n})\in Y$.

We first use induction to show that $T_{n}$ is a standard immaculate tableau, namely $T_n^{\leq c_n}$ is standard. It is clear that $T_1^{\leq c_1}$ is standard. If $T_{k}^{\leq c_{k}}=T_{k}^{<c_{k+1}}$ is standard, then, by Proposition
\ref{standard pro}, the partial tableau $T_{k+1}^{\leq c_{k+1}}=\mathrm{jdt}_{e_{k+1}}(T_{k})^{\leq c_{k+1}}$ is standard. By induction, we see that $T_n^{\leq c_n}$ is standard.

We continue to show that $S_n$ is a hook tableau. For the cell $c_{k+1}=(i,j)$ of the diagram of $\alpha$, we have $(S_n)_{ij}=(S_{k+1})_{ij}$ by the construction of $S_n$. Suppose that the path of $e_{k+1}$ in $T_{k}$ starts at $c_{k+1}=(i,j)$ and ends at $(i',j')$. If $i'>i$, then we must have $j=1$ by Lemma \ref{lemma-path}.
From \eqref{eq-sk} and the definition of the hook length it immediately follows that $(S_n)_{ij}=(S_{k+1})_{ij}\leq h_{ij}$.

Let $(P,J)=(T_n, S_n)$. The map $\phi$ is defined by $\phi(T)=(P,J)$.

Now we are able to prove our main result.

%{Given a composition $\alpha$ and a cell $c=(i,j)$ in $\alpha$, the hook
%of $c$, denoted by $H_c$, has been defined above. Now we introduce a special ``distance'' between $c$ and the cell $c'={(i',j')}$ in $H_c$, denoted by $d_H(c,c')$.
%\begin{align}
%d_H(c,c') = \left\{
%             \begin{array}{lcl}
%            j'-j, &\text{if} &i'=i; \\
%            \alpha_i+\cdots+\alpha_{i'-1}+j'-1, &\text{if} &i'>i.
%             \end{array}
%        \right.
%\end{align}}

\noindent\textbf{Proof of Theorem \ref{thm-main}.} To prove that $\psi$ is a bijection, it suffices to show that $\phi$ is the inverse map of $\psi$.

We first prove that $\phi$ is the left inverse of $\psi$, which implies the injectivity of $\psi$.
Precisely, if a pair $(P,J)\in Y$ is mapped to $T\in X$ by $\psi$, then we must have $\phi(T)=(P,J)$.
By the construction of $\psi$, there exists a sequence of pairs $\{(P_i,J_i)\}_{i=1}^n$ such that
$(P,J)=(P_1,J_1)$ and $(P_n,J_n)=(T,J_n)$, where all entries of $J_n$ are $1$.
Consider the transformation from
$(P_k,J_k)$ to $(P_{k+1},J_{k+1})$. Let $L$ denote the path determined by the cell $c_{n+1-k}$ of the diagram of $\alpha$.
Let $e_{n+1-k}$ denote the entry of $P_{k+1}$ at the cell $c_{n+1-k}$. We only need to show that
the path of $e_{n+1-k}$ in $P_{k+1}$ coincides with the path $L$. But this is clear since the partial tableaux $P_k^{\leq c_{n+1-k}}$ and $P_{k+1}^{< c_{n+1-k}}$ are standard, and $P_{k+1}$ is obtained from $P_k$ by a circular right shift of the entries on the path $L$. Therefore, we have $\mathrm{jdt}_{e_{n+1-k}}P_{{k+1}}=P_k$, as the tableau $\mathrm{jdt}_{e_{n+1-k}}P_{{k+1}}$ is obtained from $P_{k+1}$ by a circular left shift of the entries on the path $L$. This implies that the pair $(P_{k+1},J_{k+1})$ will be mapped to $(P_k,J_k)$ during the construction of $\phi(T)$. Thus, we have $\phi(T)=(P,J)$.

Next, we show that $\phi$ is the right inverse of $\psi$, which implies the surjectivity of $\psi$.
Precisely, if $T\in X$ is mapped to a pair $(P,J)\in Y$ by $\phi$, then we must have $\psi(P,J)=T$.
By the construction of $\phi$, there exists a sequence of pairs $\{(T_i,S_i)\}_{i=1}^n$ such that $(T_1,S_1)=(T,S_1)$ and $(P,J)=(T_n,S_n)$, where all entries of $S_1$ are $1$. Consider the transformation from $(T_k,S_k)$ to $(T_{k+1},S_{k+1})$. Let $e_{k+1}$ denote the entry of $T_k$ at the cell $c_{k+1}$. Note that $T_{k+1}=\mathrm{jdt}_{e_{k+1}}(T_k)$. By Lemma \ref{lemma-path}, $T_{k+1}$ is obtained from $T_k$ by a circular left shift of the entries on the path of $e_{k+1}$ in $T_k$. Moreover, the entry of $c_{k+1}$ in $S_{k+1}$ is determined by \eqref{eq-sk}, and this entry will uniquely determine a path $L$ in $T_{k+1}$ when encountering the pair $(T_{k+1},S_{k+1})$ during the construction of $\psi(P,J)$. According to the construction of the map $\psi$, the path $L$ must coincide with the path of $e_{k+1}$ in $T_k$. This implies that the pair $(T_{k+1},S_{k+1})$  will be transformed to $(T_k,S_k)$ during the construction of $\psi(P,J)$. Therefore, we have $\psi(P,J)=T$.

Combining the above two aspects, we complete the proof of the bijectivity of $\psi$.
\qed

\noindent{\bf Acknowledgements.} This work was supported by the 973 Project, the PCSIRT Project of the Ministry of Education and the National Science Foundation of China.

\end{document}